\documentclass[12pt,reqno]{amsart} 
\usepackage{amsthm,amsfonts,amssymb,amscd}

\newtheorem{theorem}[subsection]{Theorem}
\newtheorem{proposition}[subsection]{Proposition}
\newtheorem{prop}[subsection]{Proposition}

\newtheorem{question}[subsection]{Question}
\newtheorem{lemma}[subsection]{Lemma}
\newtheorem{corollary}[subsection]{Corollary}

\theoremstyle{definition}
\newtheorem{definition}[subsection]{Definition}

\newtheorem{proposition-definition}[subsection]{Proposition-Definition}

\theoremstyle{remark}
\newtheorem{remark}[subsection]{Remark}

\newcommand{\mt}[1]{\operatorname{#1}}

\newcommand{\Supp}{\operatorname{Supp}}
\newcommand{\mult}{\operatorname{mult}}

\newcommand{\GL}{{GL}_3(\CC)}
\newcommand{\SL}{{SL}_3(\CC)}
\newcommand{\PGL}{{PGL}_3(\CC)}

\newcommand{\down}[1]{\lfloor #1\rfloor}

\newcommand{\CC}{{\mathbb C}}
\newcommand{\RR}{{\mathbb R}}
\newcommand{\ZZ}{{\mathbb Z}}
\newcommand{\QQ}{{\mathbb Q}}
\newcommand{\PP}{{\mathbb P}}
\newcommand{\NN}{{\mathbb N}}

\newcommand\alp{\alpha}

\newcommand\om{\omega}
\newcommand\epsi{\epsilon}

\newcommand\lra{{\longrightarrow}}

\renewcommand\square{\frame{\phantom{{\large x}}}}

\author{D. Markushevich}

\address{D. M.: Math\'ematiques - b\^{a}t. M2, Universit\'e Lille 1,
F-59655 Villeneuve d'Ascq Cedex, France}
\email{markushe@gat.univ-lille1.fr}

\author{Yu.~G.~Prokhorov}

\thanks{The second author was partially supported
by the Russian Foundation of Fundamental Research grant
96-01-00820 and by a grant PECO-CEI from the Ministry of Higher Education
of France}

\address{Yu. G.: Algebra Section, Dept. of Mathematics,
Moscow State University, 117234 Moscow, Russia}
\email{prokhoro@mech.math.msu.su\qquad \\
prokhoro@math.pvt.msu.su\\}


\title{Exceptional quotient singularities}

\begin{document}
\maketitle

\section{Introduction and statement of main results}
The notion of exceptional singularity
was introduced by Shokurov \cite{Sh1}. A singularity $(X,P)$ is called
{\it exceptional}, if for any log canonical boundary, there is at most
one exceptional divisor of discrepancy $-1$ over $P$ 
(see Definition \ref{def-exc}). The reason for
distinguishing these singularities is that
 they have
more complicated multiple anticanonical
systems $|-nK_X|$, than
the nonexceptional ones. However, Shokurov 
suggests, and proves in dimension $3$, that
the exceptional singularities are in a sense bounded (loc. cit., Corollary 7.3).
The search of ``good" divisors in $|-nK_X|$, 
or so called $n$-{\it complements}, is an essential ingredient of 
Shokurov's project of the inductive study of log flips, log contractions
and of the classification of log canonical singularities \cite[Sect.~7]{Sh1}.

According to \cite[5.2.3, 5.6]{Sh}, \cite[1.5]{Sh1}, the exceptional
log terminal singularities in 
dimension 2 are exactly the singularities of types
$E_6,E_7,E_8$ in the sense of Brieskorn \cite{Br}. Shokurov's approach
gives a description of the dual graphs of their resolutions (cf \cite{I},
\cite[Chapter 3]{Ut}). According to \cite{Br}, the exceptional graphs correspond
to finite subgroups of $GL_2(\CC )$ of tetrahedral ($E_6$), 
octahedral ($E_7$) or icosahedral ($E_8$) types.
We found a classification-free approach to the
proof which works also in dimension 3. The exceptional
groups of types $E_6$, $E_7$, $E_8$ are exactly those finite subgroups of
$GL_2(\CC )$ which have no semiinvariants of degree $\leq 2$. Hence the 
following proposition implies Shokurov's statement.

\begin{proposition}
\label{Sh_2}
A two-dimensional quotient singularity $X=\CC^2/G$ 
by a finite group $G$ without reflections is
exceptional if and only if $G$ has no semiinvariants of degree $\le 2$.
\end{proposition}

In dimension 3, we started the study of exceptional quotient singularities
in our previous publication \cite{MP}, where we showed that the quotients
of $\CC^3$ by Klein's simple group of order 168 and by its unique 
central extension contained in $\SL$ (of order 504) are exceptional,
in using the configuration of the action of Klein's group on $\PP^2$
\cite{Kl}, \cite{W}.
The following Theorem is the main result of the present paper.

\begin{theorem}
\label{main}
A three-dimensional quotient singularity $X=\CC^3/G$ 
by a finite group $G$ without reflections is
exceptional if and only if $G$ has no semiinvariants of degree $\le 3$.
\end{theorem}

Using Miller--Blichfeldt--Dickson classification \cite{MBD}
of finite subgroups of $\GL$, we obtain a complete list of
such subgroups yielding exceptional singularities (Theorem \ref{list-exc}).

Section \ref{s-1} contains basic definitions and preliminary results. 
In particular, we
show that the exceptional divisor with discrepancy $-1$ for an exceptional
singularity is birationally unique and that its image is a point, independently
of the choice of the boundary and in any dimension. Section \ref{s-2} contains
the proofs of Proposition \ref{Sh_2}, Theorem \ref{main}, and the list of
finite subgroups of $\GL$ with exceptional quotients.
\par\noindent
{\sc Acknowledgement.}
The second author is grateful to V.~V.~Shokurov 
for  discussions on  complements and the exceptionality.

\section{Basics on exceptional singularities}
\label{s-1}

For reader's convenience, we reproduce here some basic facts
about log canonical singularities.
We are using the terminology and notation of \cite{MP}, \cite{Ut},
\cite{Sh} and \cite{Sh1} (see also \cite{KoP} for a nice introduction
to the subject).

\begin{definition}
Let $(X\ni P)$ be a normal singularity (not necessarily isolated) and 
let $D=\sum d_iD_i$ be a divisor on $X$ with real coefficients.
$D$ is called a boundary if $0\leq d_i\leq 1$ for all $i$.
It is called a subboundary, if it is majorated by a boundary.
A proper birational morphism $f \colon Y\lra X$ is called a log resolution of
$(X,D)$ at $P$, if $Y$ is nonsingular near $f^{-1}(P)$ and 
$\Supp (D)\cup E$ is a normal crossing divisor
on $Y$ near $f^{-1}(P)$, where $D$ is used
to denote both the subboundary on $X$ and its proper transform on $Y$,
and $E=\cup E_i$ is the exceptional divisor of $f$.
The pair $(X,D)$ or, by abuse of language, the divisor $K_X+D$
is called {\it terminal, canonical, Kawamata log terminal,
 purely log terminal,} and, respectively, 
{\it log canonical} near $P$, 
if the following
conditions are verified:

(i) $K_X+D$ is $\RR$-Cartier.

(ii) Let us write for any proper birational morphism $f \colon Y\lra X$
$$
K_Y\equiv f^\ast (K_X+D)+\sum a(E,X,D)E ,
$$
where $E$ runs over prime divisors on $Y$, $a(E,X,D)\in\RR$,
and $a(D_i,X,D)=-d_i$ for each component $D_i$ of $D$. Then, for some
log resolution of $(X,D)$ at $P$ and for all prime
divisors $E$ on $Y$ near $P$, we have:
$a(E,X,D)>0$ (for terminal), $a(E,X,D)\geq 0$ (for canonical), 
$a(E,X,D)> -1$ and no $d_i=1$ (for Kawamata log terminal), 
$a(E,X,D)> -1$ (for purely log terminal, without any 
restriction on the subboundary $D$), and,
respectively, $a(E,X,D)\geq -1$ (for log canonical).
\end{definition}

The coefficients $a(E,X,D)$ are called {\it discrepancies} of $f$, or of
$(X,D)$; they depend on the discrete valuations of the
function field of $X$ associated to the prime divisors $E$ and on $D$, 
but not on the choice of $f$. We will
identify prime divisors with corresponding discrete
valuations, when speaking about `divisors $E$ over $X$'
without indicating, on which birational model $E$ is realized. The conditions 
given by inequalities in part (ii) of the above definition do not depend
on the choice of a log resolution.

\begin{definition}
Let $V$ be a normal variety and let $D=\sum d_iD_i$ be 
a $\QQ$-divisor on $V$ such that $K_V+D$ is $\QQ$-Cartier.
A subvariety $W\subset V$ is said to be a {\it center of
log canonical singularities} for $(V,D)$ if there exists 
a birational morphism from a normal variety 
$g\colon \tilde{V}\to V$ and a prime divisor $E$ on 
 $\tilde{V}$ with discrepancy $a(E,V,D)\le -1$ such that $g(E)=W$.
(The case $d_i=1$, $W=D_i$ is also possible.)
The union of all centers of log canonical singularities is called the
{\it locus of log canonical singularities} and is denoted by  $LCS(V,D)$
\cite[3.14]{Sh}.
\end{definition}

The following statement is a weak form of \cite[17.4]{Ut} (in dimension $2$
it was proved earlier by Shokurov \cite{Sh}).

\begin{theorem}
\label{conn-1}
Let $X$ be a normal projective variety and let $D=\sum d_iD_i$ be 
a $\QQ$-divisor on $X$ such that $K_X+D$ is $\QQ$-Cartier.
If $-(K_X+D)$ is nef and big, then the locus of
log canonical singularities is connected.
\end{theorem}

\begin{theorem}[{\cite[6.9]{Sh}}]
\label{conn-2}
Let $X$ be a normal projective surface and let $D=\sum d_iD_i$ be 
a boundary on $X$ such that $K_X+D$ is $\QQ$-Cartier.
If $K_X+D\equiv 0$, then the locus of
log canonical singularities has one or two connected components.
\end{theorem} 

Shokurov informed us that the previous theorem is valid in any dimension
modulo Log Minimal Model Program.

\begin{definition}[{\cite[1.5]{Sh1}}]
\label{def-exc}
Let $(X\ni P)$ be a normal singularity and 
let $D=\sum d_iD_i$ be a 
boundary on $X$ such that $K_X+D$ is log canonical. The pair $(X,D)$ is 
said to be {\it exceptional} if there exists at most one exceptional divisor 
$E$ over $X$ with discrepancy $a(E,X,D)=-1$. The singularity $(X,P)$ is said to
be {\it exceptional} if $(X,D)$ is exceptional for any $D$ whenever $K_X+D$
is log canonical. 
\end{definition}

\begin{lemma}
\label{red}
Assume that there exists a reduced divisor $S=\sum S_i$ passing through $P$ 
such that $K_X+S$ is log canonical. Then $(X\ni P)$ is nonexceptional. 
\end{lemma} 

For the proof, see Lemma~1.7 in \cite{MP}. This implies, in particular, that
any three-dimensional cDV-singularity \cite{YPG}, 
and any three-dimensional terminal singularity
is nonexceptional (see \cite[1.9]{MP}).

\begin{prop}
\label{convex}
Let $(X,P)$ be an exceptional singularity. Then the divisor $E$ from
Definition \ref{def-exc} with discrepancy $a(E,X,D)=-1$ is birationally
unique. This means, that there exists a unique discrete valuation $\nu$
of the field of rational functions $k(X)$, such that for any pair
$(X,D)$, exceptional at $P$, and for any log resolution $f:\tilde{X}
\lra X$, if there is an exceptional divisor $E$ in $\tilde{X}$ with
$a(E,X,D)=-1$, then the corresponding discrete valuation $\nu_E=\nu$.
\end{prop}

\begin{proof}
Let $\{ D_1,\ldots ,D_r\}$ be a finite set of irreducible divisors on $X$,
and $f:\tilde{X}\lra X$ a log resolution of $(X,\sum_{i=1}^r D_i)$ at $P$.
We can represent the boundaries $D=\sum_i d_iD_i$ with components from
$\{ D_1,\ldots ,D_r\}$ as the unit cube $I^r\subset \RR^r$ of
vectors $(d_i)$. Then the subset $\Lambda\subset I^r$ corresponding
to log canonical pairs $(X,D)$ is given by a finite number of linear
inequalities in $d_i$'s
$$
a(E,X,D)=a(E,X,0)-\sum_id_i\mult_E(D_i)\geq -1 ,
$$
where $E$ runs over the exceptional divisors such that $f(E)\ni P$.
Let $\partial_+\Lambda$ be the closure of $\{ (d_i)\in\partial\Lambda\mid
0<d_i<1\;\forall \; i=1,\ldots ,r\}$. 
Then the exceptionality of $(X,P)$ implies that there are
no points in $\partial_+\Lambda$ satisfying the equality
$a(E,X,0)-\sum_i\mult_E(D_i)= -1$ for two different $E$'s. So,
$\partial_+\Lambda$ (if nonempty) is an open convex polyhedron lying in exactly
one hyperplane $\Pi_E$ with equation $a(E,X,0)-\sum_id_i\mult_E(D_i)= -1$.

Now, let $(X,P)$ be an exceptional singularity,
$(X,D^{(1)}),(X,D^{(2)})$ two exceptional
pairs, $f_i:X_i\lra X$ two log resolutions with
divisors $E_i$ such that $a(E_i,X,D^{(i)})=-1$. Then we can apply the previous
argument, taking the union of all the components of
$D^{(1)},D^{(2)}$ as $\{ D_1,\ldots ,D_r\}$, and some log resolution
dominating $f_1,f_2$ as $f$. 
Then we will have : firstly, by Lemma \ref{red}, 
$D^{(1)},D^{(2)}$ are represented by two different points of
$\partial_+\Lambda$,
and secondly, $\Pi_{E_1}=\Pi_{E_2}$. But then, 
both of the $E_i$ are of discrepancy
$-1$ with respect to any $D$ represented by a point of
$\partial_+\Lambda$. By exceptionality, $E_1=E_2$. This proves the
assertion of the proposition.
\end{proof}

\begin{prop}
\label{LCS-point}
Let $(X,P)$ be an exceptional singularity. Then for any pair
$(X,D)$ exceptional at $P$ and such that $LCS(X,D)\neq\varnothing$,
we have $LCS(X,D)=\{ P\}$.
\end{prop}

\begin{proof}
By Proposition \ref{convex}, the exceptional divisor $E$ with discrepancy
$-1$ is biratonally unique. Suppose, $LCS(X,D)\supsetneqq\{ P\}$.
Then for a generic hyperplane section $H$ of $(X,P)$, we have
$\mult_E(H)=0$. This implies that $a(E,X,D+\epsilon H)=a(E,X,D)=-1$, 
and, by exceptionality, $a(E',X,D+\epsilon H)>-1$ for any sufficiently small
positive $\epsi$ and any exceptional divisor $E'\neq E$. 
Let $\epsi_\ast =
\sup\;\{\epsi\;\mid\; (X,D+\epsilon H)\ \text{is log canonical}\ \}$. 
Then $(X,D+\epsi_\ast H)$
has two different divisors with discrepancy $-1$: $E$ and 
either the one of the 
$E'\neq E$ with minimal discrepancy, or $H$ in the case when $\epsi_\ast =1$.
This contradicts to the exceptionality of $(X,P)$.
\end{proof}

\begin{prop}
\label{cartesian}
Let $(X_i,P_i)$, $i=1,2,$
be two germs of normal varieties of dimension $>0$. Then 
$(X_1\times X_2,(P_1,P_2))$
is nonexceptional.
\end{prop}

\begin{proof}
If at least one of $X_i$ is not log canonical with zero boundary,
we are done, for $X_1\times X_2$ is not log canonical
with any boundary. If $(X_i,0)$ are log canonical,
then, starting from any nonzero boundaries $D^{(i)}$ at $P_i$,
we can find $\epsi_i\geq 0$ such that both pairs $(X_i,\epsi_iD^{(i)})$ 
are log canonical and possess
at least one discrete valuation of discrepancy $-1$ for $i=1,2$.
Let $W_i$ be their centers. Then 
$LCS(X_1\times X_2,\epsi_1 D^{(1)}\times X_2)$
contains $W_1\times X_2$, and
$LCS(X_1\times X_2, X_1\times \epsi_2 D^{(2)})$ 
contains $X_1\times W_2$. Hence
$(X_1\times X_2,(P_1,P_2))$
is nonexceptional by Proposition \ref{LCS-point}.
\end{proof}

\section{Quotients $\CC^m/G$}
\label{s-2}
In the first part of this section we work 
in arbitrary dimension. The following Lemma was proved in \cite[2.1]{MP} in the
$3$-dimensional case; now we generalize it to arbitrary dimension.

\begin{lemma}\label{quot-lem1}
Let $\pi\colon V\to X$ be a finite morphism of normal varieties, 
\'etale in codimension $1$. Let $D$ be a boundary on $X$, and $D'=\pi^\ast D$.
Assume that $K_X+D$ is log canonical (and then automatically is $K_V+D'$).
Then $(X,D)$ is exceptional if and only if $(V,D')$ is. 
\end{lemma}

\begin{proof}
For every proper birational morphism
$f\colon Y\to X$ consider the following commutative diagram
\begin{equation}
\label{1}
\begin{CD}
W@>{\varphi}>>Y\\
@V{g}VV@VfVV\\
V@>{\pi}>>X.\\
\end{CD}
\end{equation}
where $W$ is the normalization of the dominant component 
of $V\times_{X}Y$.
By the ramification formula (cf. \cite[proof of 3.16]{KoP})
for every exceptional divisor $F$ on $Y$ and every
exceptional divisor $E$ on $W$ which dominates $F$ we have
\begin{equation}\label{1.5}
a(E,V,D')+1=r(a(F,X,D)+1),
\end{equation}
where $r$ is the ramification index at the generic point of $E$. 
Therefore 
\begin{equation}
\label{iff}
a(E,V,D')=-1\qquad\text{ if and only if}\qquad a(F,X,D)=-1.
\end{equation}
Thus if $(X,D)$ is nonexceptional, then $(V,D')$ is so.

Conversely, assume that $(V,D')$ is nonexceptional.
Then there are infinitely many exceptional divisors over 
$V$ with discrepancy $a(\phantom{E},V,D')=-1$.
As in \cite[proof of 3.16]{KoP} we can see that 
all these divisors appear in some commutative square (\ref{1})
(for suitable $f$), whence $(X,D)$ is also nonexceptional
by (\ref{iff}).
\end{proof}

\subsubsection*{Warning}
The Lemma does not imply that the singularity of $X$ is exceptional if and
only if $V$ is. The assertion concerns {\em pairs} $(V,D')$ with
$D'$ a pullback of a boundary from $X$.

\subsection{}\label{V,G}
Let $\pi\colon V\to X$ be the quotient morphism, where $V=\CC^m$.
Let $D$ be a boundary on $X$ and let $D':=\pi^*D$. By \cite[2.2]{Sh},
(see also \cite[20.3]{Ut}, \cite[3.16]{KoP}) $K_X+D$ is log canonical 
(resp., purely log terminal, Kawamata log terminal) if and only if so 
is $K_{V}+D'$. 

\begin{question}
\label{conj}
In the notation of \ref{V,G}, assume that $G$ 
has a semiinvariant of degree
$\le m$. Does this imply that $(X\ni P)$ is nonexceptional?
\end{question}

Let $\psi$ be such a homogeneous semiinvariant of minimal 
degree $d\le m$ and let 
$D'$ be its zero locus. 
By Lemma \ref{quot-lem1} it is sufficient, for the positive
answer to 3.3, to prove that
$K_V+D'$ is log canonical. It is clear that $D'$ is a cone.
The following Proposition is an easy particular case of 
this question.

\begin{proposition} \label{impr}
Let $G$ be maximally imprimitive, that is,
let $G$ contain a normal abelian subgroup $A$ whose character subspaces
$V_\chi$ are $1$-dimensional and form one orbit under the action of $G$. 
Then $V/G$ is nonexceptional.
\end{proposition}

\begin{proof}
$G/A$ acts by permutations of the $V_\chi$, so
$G$ has a semiinvariant of the form $x_1\cdots x_m$, where
$x_i$ are coordinate linear forms. Since
$D'=\{ x_1\cdots x_m=0\}$ is a simple normal
crossing divisor, $(V,D')$ is log canonical. Then $(X,D)$ is also
log canonical, and the assertion follows
from Lemma \ref{red}.
\end{proof}

\subsection{}\label{Vi}
We cannot treat in general the case of {\it imprimitive} groups,
that is, groups $G$ which permute transitively factors
of some direct sum decomposition $V=\oplus V_i$. {\it Maximally
imprimitive} are those for which $\dim V_i =1$. But if the
action of $G$ is not transitive for some direct sum decomposition, 
the singularity is also nonexceptional. This follows from
the next proposition. 

A group $G$ is called {\it reducible}, if
$V=V_1\oplus V_2$ with $V_1,V_2$ invariant under $G$, $\dim V_i>0\; (i=1,2)$.

\begin{prop}\label{reducible}
If $G$ is reducible, then $V/G$
is nonexceptional.
\end{prop}

\begin{proof}
The proof is similar to that of Proposition \ref{cartesian},
though the result is not a corollary of \ref{cartesian}.
Let $G_i$ denote the image of $G$ in $GL(V_i)$.
Let $D_i\; (i=1,2)$ be any nonzero boundaries for $(V_i,0)$,
given by $G_i$-semiinvariant polynomials from $k[V_i]$.
We can find $\epsi_i>0$ such that both pairs $(V_i,\epsi_iD_i)$ possess
at least one discrete valuation of discrepancy $-1$ for $i=1, 2$.
Let $W_i$ be their centers. 
Define $D^{(1)}=D_1\times V_2,D^{(2)}=V_1\times D_2$ in $V$.
Then $LCS((V,\epsi_iD^{(i)})$ contains $W_1\times V_2$ or $V_1\times W_2$.
Using a commutative square of the same type as (\ref{1}), we conclude
from (\ref{1.5}) that the images of $W_1\times V_2$ and $V_1\times W_2$
in $V/G$ are also centers of log canonical singularities.
Hence $V/G$ is nonexceptional by Proposition~\ref{LCS-point}.
\end{proof}

\begin{lemma}\label{quot-lem2}
If $\dim X=m\leq 3$, the answer to \ref{conj} is affirmative.
\end{lemma}
\begin{proof} In \cite[Lemma 2.2]{MP}, we proved the assertion
in dimension $3$. The same and even easier argument
 works in dimension $2$.
\end{proof}

This Lemma shows the ``if" part of Proposition \ref{Sh_2}
and Theorem \ref{main}.

\subsection{}
\label{logic}
Now we will explain the logic of our proof of the exceptionality
of singularities in Proposition \ref{Sh_2} and Theorem \ref{main}.
Assume that $(X\ni P)$ is nonexceptional. Then
there exists a nonexceptional log canonical $K_X+D$.
Further, we will use notations of Lemma~\ref{quot-lem1}. 
Take the smallest $n\in\NN$ such that $F:=nD'$ is an 
integer divisor. Then $F$ locally near $0$ can be defined by a 
semiinvariant function, say $\psi$.
Denote $d:=\mt{mult}_0(\psi)$. 
\label{begin}
Let $\sigma\colon W\to V=\CC^m$ be the blow-up of the origin
and let $S\simeq\PP^{m-1}$ be the exceptional divisor.
Then $K_W=\sigma^*K_V+(m-1)S$ and
$\sigma^*F=R+dS$, where $R$ is the proper transform of $F$.

Further
\begin{equation}
K_W+S+\frac{m}{d}R=\sigma^*(K_V+\frac{m}{d}F).
\end{equation}
By \cite[Lemma 3.10]{KoP} $K_V+\frac{m}{d}F$ is log canonical if and only if so 
is 
$K_W+S+\frac{m}{d}R$. In this case the pair $(V,\alpha F)$ is 
exceptional for all $0\leq\alpha\leq \frac{m}{d}$ if and only if 
$K_W+S+\frac{m}{d}R$ is purely log terminal.
By the inversion of adjunction (see \cite[5.13]{Sh}, \cite[17.6]{Ut})
the purely log terminal condition for $K_W+S+\frac{m}{d}R$ is equivalent to 
that $K_S+\frac{m}{d}C$ is Kawamata log terminal, where $C=R\cap S$.
It is clear that $C$ is given by the equation 
$\psi_{{\min}}=0$, where $\psi_{{\min}}$ is 
the homogeneous component of $\psi$ of minimal degree $d$.
Therefore we have

\begin{proposition}
\label{check}
In the above notations, if $K_S+\frac{m}{d}C$
is Kawamata log terminal, then $(V,\alpha F)$ is exceptional for any
$0\leq\alpha \leq\frac{m}{d}$.
\end{proposition}\bigskip

\begin{remark}\label{remark}
In the above notations, if $G$ has no 
semiinvariants of degree $\le m$, then $\down{\frac{m}{d}C}=0$.
\end{remark}\bigskip

\subsection{Two-dimensional case}
\subsubsection*{Proof of Proposition \ref{Sh_2}}
We obtain Proposition \ref{Sh_2} as an easy corollary of \ref{check}.
Assume that $G$ has no invariants of degree $\le 2$.
Recall that $S\simeq\PP^1$ in our case, so $C$ is a finite set. 
By \ref{check} it is sufficient to prove that
$K_S+\frac{2}{d}C$ is Kawamata log terminal, and this is equivalent to that
$\down{\frac{2}{d}C}=0$. The last assertion follows by \ref{remark}.

\hspace*{\fill}\square
\medskip

\subsection{Three-dimensional case}
\subsubsection*{Proof of Theorem \ref{main}}
Assume that $G$ has no invariants of degree $\le 3$.
We shall prove that $\CC^3/G$ is exceptional.
By Proposition \ref{check} it is sufficient to prove that $K_S+\frac{3}{d}C$
is Kawamata log terminal. 
Take $c$ to be the log canonical threshold of $(S,C)$, that is, 
the maximal $\alp$ such that $K_S+\alp C$ is log canonical. 
If $K_S+\frac{3}{d}C$
is not Kawamata log terminal, then $c\le 3/d$.
First we consider the case when $c<3/d$.
Then $-(K_S+cC)$ is ample. By connectedness theorem \ref{conn-1},
the locus of log canonical singularities is connected.
By \ref{remark} $\down{\frac{3}{d}C}=0$, so any divisor of discrepancy
$-1$ should be exceptional. 
Therefore the locus of log canonical singularities on $S$ is 
a unique point, which must be invariant under the action of $G$.
The dual line gives us a semiinvariant of degree $1$, a contradiction with
our assumptions.
In the case when $c=3/d$, we can use
Theorem \ref{conn-2}. Similarly to the above, we see that
the locus of log canonical singularities is one or two points.
In both cases, there is an invariant line. This ends the
proof of Theorem \ref{main}.

\hspace*{\fill}\square
\medskip

The finite subgroups of $G\subset\GL$ were classified by
Miller--Blichfeldt--Dickson \cite{MBD} modulo central extensions
(compare with \cite{P}). There are 10 types
of such groups, denoted by A,B, \ldots ,J in \cite{MBD}.
A stands for abelian, B for reducible, and C, D are imprimitive.
The groups of type C are called {\it tetrahedral}; their image
in the symmetric group ${\frak S}_3$ permuting the 
$V_i,\; i=1,2,3$ (notations as in
\ref{Vi}) is cyclic of order 3.
The groups of type D are called {\it general monomial}; their map
to ${\frak S}_3$ is surjective. The primitive subgroups of $\GL$
belong to the $6$ types E, F, G, H, I, J.
The orders of the associated collineation groups 
$PG=G/(G\cap\CC^\ast)\subset\PGL$ 
are $36$, $72$, $216$, $60$, $360$, $168$; 
the first three are solvable,
and the last three are simple. The collineation 
groups from G to J have their names:
the Hessian group, the icosahedral one, the alternating group ${\frak A}_6$
of degree $6$, and, finally,
Klein's simple group. We have also $P$E$\triangleleft P$F$\triangleleft P$G.
\medskip

The following assertion is a consequence of Theorem \ref{main}.

\begin{theorem}\label{list-exc}
Let $G$ be a finite subgroup of $\GL$. Then the quotient $\CC^3/G$
is exceptional if and only if $G$ belongs to one of the $4$ types
{\em F, G, I, J}.
\end{theorem}

\begin{proof}
We can eliminate the types A--D by Propositions \ref{impr}, \ref{reducible}.
Further, the icosahedral group in its $3$-dimensional
representation has an invariant of degree two, which is a semiinvariant
of any group of type H. This follows, for example, from the fact that
the $3$-dimensional representation is the 
complexification of the standard real one,
which has an invariant scalar product.

\begin{lemma}\label{groupE}
Any group of type {\em E} has two nonproportional 
semiinvariants of degree $3$.
Groups of types {\em F, G} have no semiinvariants of degree $\leq 3$.
\end{lemma}

\begin{proof}
According to \cite[Sect. 115]{MBD}, a group $G$ 
of one of the types E, F, G is an extension
of a group $H$ of type D which leaves invariant the set of 
four triangles $t_1,t_2,
t_3,t_4$ defined in appropriate coordinates by the equations
$$
x_1x_2x_2=0,\; (x_1+x_2+x_3)(x_1+\om x_2+\om^{2+i}x_3)
(x_1+\om^{2}x_2+\om^{1+i}x_3)=0 \; $$
$$\phantom{x_1x_2x_2=0,\; (x_1+x_2+x_3)(x_1+\om x_2}
(\om =\exp \frac{\scriptstyle 2\pi\sqrt{-1}}
{\scriptstyle 3}, 
i=0,1,2),
$$
and the associated collineation group $PG$ is completely characterized
by its image in the group ${\frak S}_4$ of permutations on the set $\{ t_1,t_2,
t_3,t_4\}$. It is the subgroup of order $2$, conjugate to
$\{ 1, (t_1t_2)(t_3t_4)\}$ for the type E, $\{ 1, (t_1t_2)(t_3t_4),
(t_1t_3)(t_2t_4),(t_1t_4)(t_2t_3)\}$ for the type F, and the full
alternating group on four letters for the type G. Moreover, one can easily
verify that the $t_i$ belong to one pencil $\Phi$ of
plane cubics (see \cite[Lemma 4.7.6]{Sp}). So, a group of type E acts on
the pencil $\Phi$ with image $\ZZ /2\ZZ$ in $\mt{Aut}(\PP^1)$. Any involution
on $\PP^1$ has two fixed points, which implies the result for type E.

Since the image of $G$ in $\mt{Aut}(\PP^1)$ has no fixed
points for the groups of types F, G, they 
do not have semiinvariants in the pencil $\Phi$. However, any semiinvariant
of $G$ should be also that of $H$. According to \cite[Sect. 113]{MBD},
$PH$ is generated by the cycle $c=(x_1x_2x_3)$, transposition $\tau =(x_2x_3)$
and dilatation $\kappa =\:\mt{diag}(1,\om ,\om^2)$. 
Any semiinvariant of ${\frak S}_3=<c,\tau>$ is a 
polynomial in elementary symmetric
functions $\sigma_i\; (i=1,2,3)$ in $x_1,x_2,x_3$ and
in $\Delta =(x_1-x_2)(x_2-x_3)(x_1-x_3)$. A direct verification
shows that in degrees $\leq 3$, only the linear combinations of
$x_1^3+x_2^3+x_3^3=\sigma_1^3-3\sigma_1\sigma_2+3\sigma_3$
and of $\sigma_3$ are (skew-)symmetric semiinvariants 
under $\kappa$. This yields exactly the pencil $\Phi$.

\end{proof}

Thus, by Lemma \ref{quot-lem2}, the quotients of types E, H are not
exceptional. By Theorem \ref{main}, the quotients of types F, G
are exceptional. It remains to verify that the groups of types
I, J have no semiinvariants of degree $\leq 3$.

For the groups of type I, we can take a representative $I_0\subset\SL$ of
order $1080$. The homogeneous semiinvariants of any group $G$ of type I will
coincide with those of $I_0$, because they are central extensions
of the same group $PI_0\simeq {\frak A}_6$. But all the semiinvariants
of the group $I_0$ are indeed invariants. This follows from the fact that it
has no normal subgroups with abelian quotient: the only nontrivial normal
subgroup of $I_0$ is its center, isomorphic to $\ZZ /3\ZZ$, and its quotient
is the simple group ${\frak A}_6$. The algebra of invariants of $I_0$
was determined by Wiman \cite{Wi} (see also a modern account of invariants
of finite subgroups of $\SL$ in \cite{YY}, where I of \cite{MBD}
is denoted by L); it is generated by basic invariants of
degrees $6$, $12$, $30$ and $45$ with one relation of weighted 
degree $90$ between them.
Thus, there are no semiinvariants of degree $\leq 3$, and we are done.

Klein's simple group has a representation $J_{168}\subset\SL$. As above,
the only its semiinvariants are invariants, and they were determined
by Klein \cite{Kl}, see also \cite{W} or \cite{YY}. 
The degrees of basic invariants here are $4$, $6$, $14$, $21$, 
and there is one
relation of weihgted degree $42$ between them. Again there are no semiinvariants
of degree $\leq 3$, and we are done. This ends the proof of Theorem 
\ref{list-exc}.
\end{proof}

Using the classification of finite subgroups in $\SL$ (\cite{P}, \cite{YY}),
one can get the following assertion.

\begin{corollary}
Let $G$ be a finite subgroup of $\GL$ such that 
 $\CC^3/G$
is an exceptional canonical singularity. 
Then $G$ is, up to conjugation, one of the following subgroups of $\SL$:
\begin{enumerate}
\renewcommand\labelenumi{{\rm (\roman{enumi})}}
\item
Klein's simple group $J_{168}\subset\SL$,
\item
the unique central extension $J'_{504}$ of $J_{168}$ contained in 
$\SL$,
\item
the Hessian group $G_{648}\subset\SL$,
\item
the normal subgroup $F_{216}$ of $G_{648}$,
\item
a central extension $I_{1080}$ of $\mathfrak{A}_6$.
\end{enumerate}
\end{corollary}

\begin{proof}
It suffices to show that $G\subset\SL$. Let $r$ be the order of
the center $Z$ of $G$. Then $\CC^3/Z$ is resolved by a single
blow up, giving an exceptional divisor with discrepancy
$-1+\frac{3}{r}$. Construct, as in the proof of Lemma \ref{quot-lem1},
a diagram (\ref{1}) with $\CC^3/Z,\: \CC^3/G$ in place of $V,X$ respectively.
Then the formula (\ref{1.5}) implies that the minimal discrepancy
of exceptional divisors $E$ over $X\;$ $\;\; a_{\min}\leq -1+\frac{3}{r}$.
For canonical singularities, $a_{\min}\geq 0$, hence $r\leq 3$. If
$r=1$ or $3$, $\; G\subset\SL$ and we are done.

Assume that $r=2$. Then $Z=\{ \pm 1\}$. Let $G_0=G\cap\SL$. 
The orders of the subgroups of $\SL$ of types
F, G, I, J are even, hence $G_0$ contains an element $g$ of 
order 2. Then either $g$ or $-g$ is a reflection.
This contradicts to our hypotheses. Hence $r=2$ is impossible.
\end{proof}

Recall that a variety $X$ is said to have 
$\epsilon$-log terminal singularities if 
$a(E,X,0)>-1+\epsilon$ for any exceptional divisor $E$ over $X$.
Adapting the proof of the previous corollary to the case
$a_{\min} > -1+\epsi$ instead of $a_{\min}\geq 0$, we obtain the
following result:

\begin{corollary}
Fix $\epsilon>0$. Then the set of subgroups $G\subset\GL$
without reflections such that $\CC^3/G$ is an exceptional 
$\epsilon$-log terminal singularity is finite up to 
conjugation.
\end{corollary}

A similar finiteness result was obtained by \cite{Bor} for abelian quotients
of any dimension (which are never exceptional by Proposition \ref{reducible}).

\end{document}